 \documentclass[12pt,reqno]{article}
\usepackage{amsfonts,amssymb,eucal,amsmath,amssymb}

\date{}

\parfillskip=\parindent plus1fil %prevents last line of paragraph being full

\renewcommand{\phi}{\varphi}

\renewcommand{\emph}{\textsl}

\newcommand{\qed}{\parfillskip=0pt plus1fil\nolinebreak\hfill$\rule{1ex}{1ex}$%
\par\addvspace{12pt plus 3pt minus 3pt}\parfillskip=\parindent plus1fil}

\title{\sffamily\bfseries\Large Essential norm of Toeplitz operators and  Hankel operators on the weighted Bergman space}
\author{\sffamily\bfseries Fengying Li\footnote{Supported by China Scholarship Council}}

\newenvironment{keeptogether}{\pagebreak[0]\samepage}{}

\numberwithin{equation}{section}

\newtheorem{Theorem}[equation]{\sffamily\bfseries Theorem}
\newtheorem{Proposition}[equation]{\sffamily\bfseries Proposition}
\newtheorem{Lemma}[equation]{\sffamily\bfseries Lemma}
\newtheorem{Corollary}[equation]{\sffamily\bfseries Corollary}
\newtheorem{Conjecture}[equation]{\sffamily\bfseries Conjecture}

\newenvironment{theorem}{\begin{keeptogether}\begin{Theorem}\slshape}%
  {\end{Theorem}\end{keeptogether}}
  {\end{Proposition}\end{keeptogether}}
  {\end{Lemma}\end{keeptogether}}
  {\end{Corollary}\end{keeptogether}}
  {\end{Conjecture}\end{keeptogether}}

\newcounter{subtheoremc}
\newcounter{subsubtheoremc}

\newenvironment{subtheorem}{\begin{list}{{\normalfont
(\alph{subtheoremc})}\hfill}{\usecounter{subtheoremc}
\setlength{\labelwidth}{.4in}
\setlength{\labelsep}{0pt}
\setlength{\leftmargin}{.4in}
}}{\end{list}}

\makeatletter

\renewcommand{\@startsection}[6]{\if@noskipsec \leavevmode
\fi
   \par \@tempskipa #4\relax
   \@afterindenttrue
   \ifdim \@tempskipa <\z@ \@tempskipa -\@tempskipa \@afterindenttrue\fi
   \if@nobreak \everypar{}\else
     \addpenalty{\@secpenalty}\addvspace{\@tempskipa}\fi \@ifstar
     {\@ssect{#3}{#4}{#5}{#6}}{\@dblarg{\@sect{#1}{#2}{#3}{#4}{#5}{#6}}}}

\renewcommand{\section}{\@startsection {section}{1}{\z@}%
                                       {-3.5ex \@plus -1ex \@minus -.2ex}%
                                       {2.3ex \@plus.2ex}%
                                       {\reset@font\large\sffamily\bfseries}}

\makeatother

\renewcommand{\[}{\parfillskip=0pt plus1fil$$}
\renewcommand{\]}{$$\parfillskip=\parindent plus1fil}

\makeatletter
\renewenvironment{equation}{\parfillskip=0pt plus1fil%
\refstepcounter{equation}%
$$%
}{%
\leqno\tagform@{\theequation}%
$$%
\parfillskip=\parindent plus1fil}%

\newcommand{\simpletag}[1]{\def\@currentlabel{#1}\def\theequation{#1}%
  \addtocounter{equation}{-1}}
\makeatother

\begin{document}

\maketitle \thispagestyle{empty}

 \textbf{Abstract.} In this paper,
we show that  on the weighted Bergman space of the unit disk the
essential norm of a noncompact Hankel operator equals its distance
to the set of compact Hankel operators and is realized by infinitely
many compact Hankel operators, which is analogous to the theorem of
Axler, Berg, Jewell and Shields on the Hardy space in \cite{AxD};
moreover, the distance is realized by infinitely many compact Hankel
operators with symbols continuous on the closure of the unit disk
and vanishing on the unit circle.



\section{Introduction}
{}Let $D$ be the open unit disk in the complex plane $\mathbb{C}$.
    Let $L^{\infty}(D)$ denote the space of bounded measurable functions on the unit disk $D$, and let $H^{\infty}(D)$ denote its subalgebra of bounded analytic functions. We write $dA$ to denote the normalized Lebesgue area measure on the unit disk $D$. For $\alpha>-1$, $L^{2}(D,dA_{\alpha})$ consists of all function $f$ on $D$ such that
     $$\|f\|_{2,\alpha}^{2}=\int_{D}|f(z)|^{2}dA_{\alpha}(z)<\infty$$ where $dA_{\alpha}(z)=(1+\alpha)(1-|z|^{2})^{\alpha}dA(z)$, it easy to see that $\int_{D}dA_{\alpha}(z)=1$.

    For $\alpha>-1$ we define the \textit{weighted Bergman space}
    $$A_{\alpha}^{2}(D)=H(D)\cap L^{2}(D,dA_{\alpha})$$
    where $H(D)$ is the space of analytic functions on $D$.  When $\alpha=0$ we recover the standard definition of the Bergman space.

      The \textit{Toeplitz operator} $T_{\varphi}$ with symbol $\varphi \in L^{\infty}(D)$ on the weighted Bergman spaces is defined by
            $$T_{\varphi}f(z)=P_{\alpha}(\varphi f)(z)=\int_{D}K_{\alpha}(z,w)f(w)\phi(w)dA_{\alpha}(w) \ \ \forall f\in A_{\alpha}^{2}(D)$$
             where $K_{\alpha}(z,w)=\frac{1}{(1-z\overline{w})^{2+\alpha}}$ is the \textit{reproducing kernel} and $P_{\alpha}$ is the orthogonal projection of $L^{2}(D,dA_{\alpha})$ onto $A^{2}_{\alpha}(D)$. Theorem 4.24 in \cite{Zhu} guarantees $P_{\alpha}$ is bounded operator on the weighted Bergman space.

             The \textit{Hankel operator} $H_{\varphi}$ with symbol $\varphi \in L^{\infty}(D)$ on the weighted Bergman space is defined by
    $$H_{\varphi}f=(I-P_{\alpha})(\varphi f) \ \ \forall f\in A^{2}_{\alpha}(D)$$
    where $I-P_{\alpha}$ is the orthogonal projection from $L^{2}(D,dA_{\alpha})$ onto $(A_{\alpha}^{2}(D))^{\perp}$.

     Let $\mathcal{K}(D)$ denote the space of compact operators on $D$. The \textit{essential norm} of an operator $T$ is defined by
           \[
           \|T\|_{e}=inf\{\|T-K\|:K\in \mathcal{K}(D)\};
           \]
           i.e., the distance $T$ to the space of compact operators.

           In order to address an approximation problem in $L^{\infty}(\partial D)$ where $\partial D$ is the boundary of unit disk $D$ (\cite{AAK}, \cite{Sar}), Axler, Berg, Jewell and Shields \cite{AxD} obtained  the following  beautiful result for the Hardy space $H^{2}(\partial D)$.
 \begin{theorem}\label{13}
  Let $H_f$ be a noncompact Hankel operator on Hardy space $H^{2}(\partial D)$.  Then there exist infinitely many different compact Hankel operators $H_\phi$ such that $\|H_f - H_\phi\|=\|H_f\|_{e}$.
 \end{theorem}

    In other words, for a noncompact Hankel operator  $H_{f}$ with symbol $f\in L^{\infty}(\partial D)$ on Hardy space $H^{2}(\partial D)$, its distance to the space of compact operators is realized by infinitely many compact  Hankel operators. On the Hardy space, a theorem of Nehari \cite{SCP}  states
    $$\|H_f\|=dist (f, H^{\infty}).$$
    It is also true that $\|H_f\|_e=dist (f, H^{\infty}+C(\partial D))$.

    For Hankel (Toeplitz) operators on unweighted Bergman space, it is known the essential norm is realized by some compact operators(\cite{GoK},\cite{HS}), and in \cite{Zhe}, the essential norm was estimated by
    \[
 \overline{\lim_{|z|\rightarrow1}}\|f\circ\varphi_{z}-P(f\circ\varphi_{z})\|_{2}\leq\|H_{f}\|_{e}\leq C\overline{\lim_{|z|\rightarrow1}}\|f\circ\varphi_{z}-P(f\circ\varphi_{z})\|_{2}^{1/10}
\]
where $\varphi_{z}(w)=\frac{z-w}{1-\overline{z}w}$ and $C$ is a constant. The estimation of the Teoplitz operator $T_{f}$ is
 \[
 \overline{\lim_{|z|\rightarrow1}}\|P(f\circ\varphi_{z})\|_{2}\leq\|T_{f}\|_{e}\leq C\overline{\lim_{|z|\rightarrow1}}\|P(f\circ\varphi_{z})\|_{2}^{1/10}.
\]

    A natural question is  whether Theorem 1.1 can be extended to the case of Hankel (Toeplitz) operators on the weighted Bergman space.

Indeed, we shall prove that the conclusion of Theorem~\ref{13} holds for noncompact Hankel operators $H_{f}$ and noncompact Toeplitz operators $T_{f}$, with symbol $f\in L^{\infty}(D)$; moreover, the symbols $\varphi$ of the approximation operators $H_{\varphi}$ ($T_{\varphi}$) will actually reside in a better behaved function space, which will guarantee that $H_{\varphi}$,
$H_{\overline{\varphi}}$, $T_{\varphi}$, and $T_{\overline{\varphi}}$ are compact.

Let $C_\partial (\overline{D})$ denote the space of continuous functions on the closure $\overline{D}$ of the unit disk and vanishing on the unit circle $\partial D.$  By results in \cite{Axl}, \cite{McS}, \cite{Str}, \cite{KZd}, \cite{Zhe}, \cite{Zhu}, \cite{Zkh}, we easily see that for each $\varphi$ in $C_\partial (\overline{D})$, $H_{\varphi}$,
$H_{\overline{\varphi}}$, $T_{\varphi}$, and $T_{\overline{\varphi}}$ are compact. The first two theorems are inspired by and analogous to Theorem \ref{13}.

\begin{theorem} \label{39}
   Let $f\in L^{\infty}(D)$, and $H_{f}$  the associated noncompact Hankel operator on $A^{2}_{\alpha}(D)$. There exist infinitely many distinct compact Hankel operators $H_{\varphi}$ with symbol $\varphi$ in $C_\partial (\overline{D})$ such that
   $$\|H_{f}-H_{\varphi}\|=\|H_{f}\|_{e}.$$
\end{theorem}

\begin{theorem} \label{14}
   Let $f\in L^{\infty}(D)$, and $T_{f}$  the associated noncompact Teoplitz operator on $A^{2}_{\alpha}(D)$. There exist infinitely many distinct compact Toeplitz operators $T_{\varphi}$ with symbol $\varphi$ in $C_\partial (\overline{D})$ such that $\|T_{f}-T_{\varphi}\|=\|T_{f}\|_{e}$.
\end{theorem}

If $f$ is harmonic on the unit disk, we have the following result.

\begin{theorem} \label{18}
   Let $f$ be a bounded harmonic function on the unit disk, and $H_{f}$  the associated noncompact Hankel operator on $A^{2}_{a}(D)$. There exist infinitely many distinct harmonic functions $\varphi$ on the unit disk and continuous on the closure of the unit disk such that $$\|H_{f}-H_{\varphi}\|=\|H_{f}\|_{e}.$$
\end{theorem}

For the Hardy space, the reproducing kernel is given by $K_{z}(\xi)=\frac{1}{1-z\overline{\xi}}$ for $z\in D$ and $\xi\in \partial D$ (see Corollary 2.11 in \cite{CM}); i.e., for any $f\in H^{2}(\partial D)$, $f(z)=<f,K_{z}>$. Since the reproducing kernel for $A_{\alpha}^{2}(D)$ is $K_{\alpha}(z,w)=\frac{1}{(1-z\overline{w})^{2+\alpha}}$ with $\alpha > -1$, it is tempting to consider the Hardy space as a limiting situation for the weighted Bergman spaces as $\alpha \rightarrow -1$; for this reason, the Hardy space is often denoted as $A^{2}_{-1}$. As a result, we view the theorems in this paper as a generalization in which the theorem of Axler, Berg, Jewell and Shields's appears as a special limiting case.

Recall that a sequence  $\{A_{n}\}$ of bounded linear operators on a Banach space $H$ converges to an operator $A$ in the \textit{strong operator topology} if $\|(A_{n}-A)f\|\rightarrow 0$ for every $f\in H$.

Our proof will make use of the following result which we will record here, and whose proof can be found in \cite{AxD}.
\begin{theorem}\label{40}
 Let $H_{1}$ and $H_{2}$ be two Hilbert spaces, and $T:H_{1}\rightarrow H_{2}$ a noncompact bounded operator. Let $\{T_{n}\}_{n\geq1}$ be a sequence of compact operators from $H_{1}$ to $H_{2}$ such that $T_{n}\rightarrow T$ and $T_{n}^{*}\rightarrow T^{*}$ in the strong operator topology. Then there exist sequences $\{a_{n}\}_{n\geq1}$ and $\{b_{n}\}_{n\geq1}$ of non-negative real numbers such that $\sum_{n\geq1}a_{n}=\sum_{n\geq1}b_{n}=1$ and $\|T-K_{1}\|=\|T-K_{2}\|=\|T\|_{e}$, where $K_{1}=\sum_{n\geq1}a_{n}T_{n}$ and $K_{2}=\sum_{n\geq1}b_{n}T_{n}$; moreover, $K_{1}\neq K_{2}$.
\end{theorem}

\section{Proof of Theorems}

In order to use Theorem~\ref{40} to prove Theorems~\ref{39} and~\ref{14}, we need to establish all the conditions in the premise; i.e., there  exists a sequence of functions $\psi_{n}\in C_{\partial} (\overline{D})$ such that the sequence of compact Hankel(Toeplitz) operators $\{H_{\psi_{n}}\}(\{T_{\psi_{n}}\})$ and
 $\{H_{\psi_{n}}^{*}\}(\{T_{\psi_{n}}^{*}\})$ respectively converge to $H_{f}(T_{f})$ and $H_{f}^{*}(T_{f}^{*})$ in the strong operator topology. We start by approximating $f\in L ^{\infty}(D)$ by continuous functions.

 Suppose $\delta$ is a positive smooth function on the complex plane $\mathbb{C}$ such that\\
  (a) $\delta$ is compactly supported and identically zero outside of $D$,\\
  (b) $\int_{\mathbb{C}}\delta(z)dA_{\alpha}(z)=1$,\\
  (c) For $\varepsilon>0$, $\lim_{\varepsilon\rightarrow 0}\delta_{\varepsilon}(z)$ is a Dirac delta function where $\delta_{\varepsilon}(z)=\frac{1}{\epsilon^{2}}\delta(\frac{z}{\epsilon})$,\\
  (d) $\int_{|z|>\varepsilon}\delta_{\varepsilon}(z)dA_{\alpha}(z)=0$.\\
  Then $\delta$ is called a \emph{mollifier} and $\int_{\mathbb{C}}\delta_{\varepsilon}(z)dA_{\alpha}(z)=1$.

  Any $f\in L^{\infty}(D)$ can be extended to the whole complex plane $\mathbb{C}$ by taking it to be zero outside of $D$. For convenience, we will denote it by the same function and so we can assume $f\in L^{\infty}_{loc}(\mathbb{C})$; thus, $f\in L^{1}_{loc}(\mathbb{C},dA_{\alpha})$ since $$\int_{\mathbb{C}}|f(z)|dA_{\alpha}(z)=\int_{D}|f(z)|dA_{\alpha}(z)\leq\|f\|_{\infty}<\infty.$$

 We can define the convolution
$$\delta_{\varepsilon}\ast f(z)=\int_{\mathbb{C}}\delta_{\varepsilon}(z-w)f(w)dA_{\alpha}(w)=\int_{\mathbb{C}}\delta_{\varepsilon}(w)f(z-w)dA_{\alpha}(w).$$
 For each fixed $z\in D$, the non-trivial domain of integration for $\int_{\mathbb{C}}\delta_{\varepsilon}(z-w)f(w)dA_{\alpha}(w)$ is the disk centered at $z$ and of radius $\epsilon$. Note the convolution is still defined for $z\in \partial D$;  hence $\delta_{\varepsilon}\ast f$ is a mollification of $f$.

 It is well known that \\
 (a) $\delta_{\varepsilon}\ast f\in C^{\infty}(\mathbb{C},dA_{\alpha})$, \\
 (b) $\delta_{\varepsilon}\ast f\in L^{2}(\mathbb{C},dA_{\alpha})$ and $\|\delta_{\varepsilon}\ast f-f\|_{2,\alpha}\rightarrow 0$ as $\varepsilon\rightarrow 0$.\\
 The reader may wish to consult \cite{Con}, \cite{DiB} and \cite{MoT} for more information.

Note, even if the function $f\in L^{\infty}(D)$ is zero on the boundary $\partial D$, the convolution $\delta_{\varepsilon}\ast f$ may not be identically zero on $\partial D$. We will need to modify the convolution to make sure that does not happen.

For a sequence $\{r_{n}\}$ such that $0<r_{n}<1$ and $z\in D$, we define
\begin{equation}\label{15}
 f_{r_{n}}(z)=\begin{cases}f(z) \ \ \ &
|z|<r_{n},\\
0 \ \ \ &
|z|\geq r_{n}.
  \end{cases}
\end{equation}

We claim that for any $\varepsilon>0$ the convolution $\delta_{\varepsilon}\ast f_{r_{n}}$ has the following properties:
\begin{subtheorem}\label{12}
\item
 $\delta_{\varepsilon}\ast f_{r_{n}}(z)$ converges to $\delta_{\varepsilon}\ast f(z)$ pointwise as $r_{n}\rightarrow 1$;

 \item
 If $dist(r_{n}D,\partial D)>\varepsilon$, the convolution $\delta_{\varepsilon}\ast f_{r_{n}}$ is equal to zero on $\partial D$.
 \end{subtheorem}

Proof of this claim:

(a) To see this, for every $z\in D $
\begin{equation}
  \begin{aligned}
  &|\delta_{\varepsilon}\ast f(z)-\delta_{\varepsilon}\ast f_{r_{n}}(z)|\\
  &\leq\int_{\mathbb{C}}\delta_{\varepsilon}(z-w)|f(w)-f_{r_{n}}(w)|dA_{\alpha}(w)\\
  &=\int_{|z-w|\leq\epsilon}\delta_{\varepsilon}(z-w)|f(w)-f_{r_{n}}(w)|dA_{\alpha}(w)\\
  &+2\|f\|_{\infty}\int_{|z-w|>\epsilon}\delta_{\varepsilon}(z-w)dA_{\alpha}(w)
  \end{aligned}
\end{equation}
Since $f_{r_{n}}(z)\rightarrow f(z)$ pointwise as $r_{n}\rightarrow 1$ and $\|f_{r_{n}}\|_{\infty}\leq\|f\|_{\infty}$, for any $z\in \mathbb{C}$ we have $\delta_{\varepsilon}(z-w)f_{r_{n}}(w)\rightarrow \delta_{\varepsilon}(z-w)f(w)$ pointwise as $r_{n}\rightarrow 1$ and $|\delta_{\varepsilon}(z-w)f_{r_{n}}(w)|\leq \delta_{\varepsilon}(z-w)\|f\|_{\infty}$. By the dominated convergence theorem, the last equality of (2.2) goes to 0 as $r_{n}\rightarrow 1$.

(b) We'll show that $\delta_{\varepsilon}\ast f_{r_{n}}|\partial D=0$. In fact,
\begin{equation}
  \begin{aligned}
&\delta_{\varepsilon}\ast f_{r_{n}}(z)|\partial D\\
&=\int_{\mathbb{C}}\delta_{\varepsilon}(z-w)f_{r_{n}}(w)dA(w)|\partial D\\
&=\int_{|z-w|<\epsilon}\delta_{\varepsilon}(z-w)f(w)\chi_{r_{n}D}(w)dA(w)|\partial D
   \end{aligned}
\end{equation}
where $\chi_{A}$ is the characteristic function of the set $A$.
 For $z\in\partial D$, from the assumptions $dist(r_{n}D,\partial D)>\epsilon$ and $|z-w|<\epsilon$, the domain of the last integration is empty. Thus $\delta_{\varepsilon}\ast f_{r_{n}}(z)|\partial D=0$ for all $0<r_{n}<1$. Hence $\delta_{\varepsilon}\ast f_{r_{n}}\in C_{\partial}(\overline{D})$. This finishes the claim.

For any $f\in L^{\infty}(D)$, we have a sequence of functions $\delta_{\varepsilon}\ast f_{r_{n}}(z)\in C_{\partial} (\overline{D})$. In \cite{Zhe} and \cite{Zhu}, it was established that $H_{\delta_{\varepsilon}\ast f_{r_{n}}}$($H_{\delta_{\varepsilon}\ast f_{r_{n}}}^{*}$) and $T_{\delta_{\varepsilon}\ast f_{r_{n}}}$ ($T_{\delta_{\varepsilon}\ast f_{r_{n}}}^{*}$) are compact on Bergman space for $0<r_{n}<1$. The next step is to prove convergence in the strong operator topology of these sequences of operators on the weighted Bergman space.\\

\textbf{Proof of theorem~\ref{39}}. First, we show that $H_{\delta_{\varepsilon}\ast f_{r_{n}}}$ converges to $H_{f}$ in the strong operator topology.

 It is well known that the subalgebra $H^{\infty}(D)$ is dense in $A_{\alpha}^{2}(D)$, i.e. $\forall \varepsilon_{1}>0$ and for any $g\in A_{\alpha}^{2}(D)$, there exists a $g_{1}\in H^{\infty}(D)$ such that $\|g-g_{1}\|_{2,\alpha}<\varepsilon_{1}$. Using the Holder inequality, the bounds $\|\delta_{\varepsilon}\ast f_{r_{n}}\|_{2,\alpha}\leq\|f\|_{\infty}$, and the fact that the orthogonal projection $I-P_{\alpha}$ is a bounded operator on the weighted Bergman space $A^{2}_{\alpha}(D)$, we see that
 \begin{equation}
  \begin{aligned}
&\|(H_{f}-H_{\delta_{\varepsilon}\ast f_{r_{n}}})g\|_{2,\alpha} \\
&\leq\|(H_{f}-H_{\delta_{\varepsilon}\ast f_{r_{n}}})(g-g_{1})\|_{2,\alpha}+\|(H_{f}-H_{\delta_{\varepsilon}\ast f_{r_{n}}})g_{1}\|_{2,\alpha}\\
&= \|(I-P_{\alpha})(f-\delta_{\varepsilon}\ast f_{r_{n}})(g-g_{1})\|_{2,\alpha}+\|(I-P_{\alpha})(f-\delta_{\varepsilon}\ast f_{r_{n}})g_{1}\|_{2,\alpha}\\
&\leq \|(f-\delta_{\varepsilon}\ast f_{r_{n}})(g-g_{1})\|_{2,\alpha}+\|(f-\delta_{\varepsilon}\ast f_{r_{n}})g_{1}\|_{2,\alpha}\\
&\leq \|f-\delta_{\varepsilon}\ast f_{r_{n}}\|_{2,\alpha}\|g-g_{1}\|_{2,\alpha}+\|(f-\delta_{\varepsilon}\ast f)g_{1}\|_{2,\alpha}\\
&+\|(\delta_{\varepsilon}\ast f-\delta_{\varepsilon}\ast f_{r_{n}})g_{1}\|_{2,\alpha}\\
&\leq 2\varepsilon_{1}\|f\|_{\infty}+\|f-\delta_{\varepsilon}\ast f\|_{2,\alpha}\|g_{1}\|_{\infty}+\|(\delta_{\varepsilon}\ast f-\delta_{\varepsilon}\ast f_{r_{n}})g_{1}\|_{2,\alpha}
 \end{aligned}
\end{equation}
From claim (a), we see that $\delta_{\varepsilon}\ast f_{r_{n}}(z)\rightarrow\delta_{\varepsilon}\ast f(z)$ pointwise as $r_{n}\rightarrow 1$.
   For $g_{1}\in H^{\infty}(D)$, we also have $\delta_{\varepsilon}\ast f_{r_{n}}(z)g_{1}(z)\rightarrow\delta_{\varepsilon}\ast f(z)g_{1}(z)$ pointwise. It is easy to see that $|\delta_{\varepsilon}\ast f_{r_{n}}(z)g_{1}(z)|\leq\|f\|_{\infty}|g_{1}(z)|$. By the dominated convergence theorem $\|(\delta_{\varepsilon}\ast f-\delta_{\varepsilon}\ast f_{r_{n}})g_{1}\|_{2,\alpha}\rightarrow 0$ as $r_{n}\rightarrow 1$; thus,
$\|(H_{f}-H_{\delta_{\varepsilon}\ast f_{r_{n}}})g\|_{2}\rightarrow 0$ as $r_{n}\rightarrow1$ and $\varepsilon\rightarrow0$ for all $\varepsilon_{1}>0$. We have shown that $H_{\delta_{\varepsilon}\ast f_{r_{n}}}$ converges to $H_{f}$ in the strong operator topology.

Next we show that $H_{\delta_{\varepsilon}\ast f_{r_{n}}}^{*}$ converges to $H_{f}^{*}$ in the strong operator topology.

For any $g\in A_{\alpha}^{2}(D)$ and $\forall\varepsilon_{1}>0$, there exists a $g_{1}\in H^{\infty}(D)$ such that $\|g-g_{1}\|_{2,\alpha}<\varepsilon_{1}$. Similar to the previous argument,
\begin{equation}\label{16}
  \begin{aligned}
   &\|(H_{f}^{*}-H_{\delta_{\varepsilon}\ast f_{r_{n}}}^{*})g\|_{2,\alpha} \\
   &\leq \|(H_{f}^{*}-H_{\delta_{\varepsilon}\ast f_{r_{n}}}^{*})(g-g_{1})\|_{2,\alpha}+\|(H_{f}^{*}-H_{\delta_{\varepsilon}\ast f_{r_{n}}}^{*})g_{1}\|_{2,\alpha}\\
   &=\|(I-P_{\alpha})^{*}(f-\delta_{\varepsilon}\ast f_{r_{n}})^{*}(g-g_{1})\|_{2,\alpha}\\
   &+\|(I-P_{\alpha})^{*}(f-\delta_{\varepsilon}\ast f_{r_{n}})^{*}g_{1}\|_{,\alpha2}\\
   &\leq \|(\overline{f-\delta_{\varepsilon}\ast f_{r_{n}}})(g-g_{1})\|_{2,\alpha}+\|(\overline{f-\delta_{\varepsilon}\ast f_{r_{n}}})g_{1}\|_{2,\alpha}\\
   &\leq \|\overline{f-\delta_{\varepsilon}\ast f_{r_{n}}}\|_{2,\alpha}\|g-g_{1}\|_{2,\alpha}+\|(\overline{f-\delta_{\varepsilon}\ast f})g_{1}\|_{2,\alpha}\\
   &+\|\overline{(\delta_{\varepsilon}\ast f-\delta_{\varepsilon}\ast f_{r_{n}}})g_{1}\|_{2,\alpha}\\
   &\leq 2\varepsilon_{1}\|f\|_{\infty}+\|(\overline{f-\delta_{\varepsilon}\ast f})\|_{2,\alpha}\|g_{1}\|_{\infty}+\|(\overline{\delta_{\varepsilon}\ast f-\delta_{\varepsilon}\ast f_{r_{n}}})g_{1}\|_{2,\alpha}
  \end{aligned}
\end{equation}
For $g_{1}\in H^{\infty}(D)$, we have $\overline{\delta_{\varepsilon}\ast f_{r_{n}}}(z)g_{1}(z)\rightarrow\overline{\delta_{\varepsilon}\ast f}(z)g_{1}(z)$ pointwise  and $|\overline{\delta_{\varepsilon}\ast f_{r_{n}}}(z)g_{1}(z)|\leq\|f\|_{\infty}|g_{1}(z)|$. By the dominated convergence theorem, we have $\|(\overline{f-\delta_{\varepsilon}\ast f_{r_{n}}})g_{1}\|_{2,\alpha}\rightarrow 0$ as $r_{n}\rightarrow 1$. Using the fact that $\|\overline{f-\delta_{\varepsilon}\ast f}\|_{2,\alpha}\rightarrow 0$ as $\varepsilon\rightarrow 0$, and for all $\varepsilon_{1}>0$, we see that the last line in (\ref{16}) goes to 0.

   From Theorem~\ref{40}, there exist two sequences $\{a_{n}\}_{n\geq1}$ and $\{b_{n}\}_{n\geq1}$ of non-negative real numbers such that $\sum_{n\geq1}a_{n}=\sum_{n\geq1}b_{n}=1$. Let $\varphi_{1}=\sum^{\infty}_{n=1}a_{n}\psi_{n}$ and  $\varphi_{2}=\sum^{\infty}_{n=1}b_{n}\psi_{n}$ where $\psi_{n}=\delta_{\varepsilon}\ast f_{r_{n}}$. Since each $\psi_{n}$ is continuous on $\overline{D}$ and $\psi_{n}|_{\partial D}=0$, $\varphi_{1}$ and $\varphi_{2}$ are continuous on $\overline{D}$ and equal to zero on $\partial D$. This implies $H_{\psi_{n}}$ is compact for $n\geq1$.
From the formula,
$$H_{\varphi_{1}}=H_{\sum^{\infty}_{n=1}a_{n}\psi_{n}}=\sum^{\infty}_{n=1}a_{n}H_{\psi_{n}}$$
 it follows that $H_{\varphi_{1}}$ is compact; similarly for $H_{\varphi_{2}}$. Theorem~\ref{40} guarantees that $H_{\varphi_{1}}\neq H_{\varphi_{2}}$.

 The two distinct compact Hankel operators $H_{\varphi_{1}}$ and $H_{\varphi_{2}}$ satisfy $\|H_{f}-H_{\varphi_{1}}\|=\|H_{f}-H_{\varphi_{2}}\|=\|H_{f}\|_{e}$. Let $\varphi=s\varphi_{1}+(1-s)\varphi_{2}$, for $s\in (0,1)$. Hence there exist infinitely many compact Hankel operators $H_{\varphi}$ such that $\|H_{f}-H_{\varphi}\|=\|H_{f}\|_{e}$. This finishes the proof of Theorem~\ref{39}.\qed

  \textbf{ Proof of Theorem~\ref{14}}. We show $T_{\delta_{\varepsilon}\ast f_{r_{n}}}$ converges to $T_{f}$ in the strong operator topology.

For any $g\in A^{2}_{\alpha}(D)$, and $\varepsilon_{1}>0$, there exists a $g_{1}\in H^{\infty}(D)$ such that $\|g-g_{1}\|_{2,\alpha}<\varepsilon_{1}$. Then using the Holder inequality, $\|\delta_{\varepsilon}\ast f_{r_{n}}\|_{2,\alpha}\leq\|f\|_{\infty}$ and the boundedness of the orthogonal projection $P_{\alpha}$ on the weighted Bergman space, we have
\begin{equation}
  \begin{aligned}
&\|(T_{f}-T_{\delta_{\varepsilon}\ast f_{r_{n}}})g\|_{2,\alpha} \\
&\leq\|(T_{f}-T_{\delta_{\varepsilon}\ast f_{r_{n}}})(g-g_{1})\|_{2,\alpha}+\|(T_{f}-T_{\delta_{\varepsilon}\ast f_{r_{n}}})g_{1}\|_{2,\alpha}\\
&=\|P_{\alpha}(f-\delta_{\varepsilon}\ast f_{r_{n}})(g-g_{1})\|_{2,\alpha}+\|P_{\alpha}(f-\delta_{\varepsilon}\ast f_{r_{n}})g_{1}\|_{2,\alpha}\\
&\leq \|f-\delta_{\varepsilon}\ast f_{r_{n}}\|_{2,\alpha}\|(g-g_{1})\|_{2,\alpha}+\|(f-\delta_{\varepsilon}\ast f_{r_{n}})g_{1}\|_{2,\alpha}\\
&\leq 2\varepsilon_{1}\|f\|_{\infty}+\|(f-\delta_{\varepsilon}\ast f)g_{1}\|_{2,\alpha}+\|(\delta_{\varepsilon}\ast f-\delta_{\varepsilon}\ast f_{r_{n}})g_{1}\|_{2,\alpha}\\
&\leq 2\varepsilon_{1}\|f\|_{\infty}+\|(f-\delta_{\varepsilon}\ast f)\|_{2,\alpha}\|g_{1}\|_{\infty}+\|(\delta_{\varepsilon}\ast f-\delta_{\varepsilon}\ast f_{r_{n}})g_{1}\|_{2,\alpha}
 \end{aligned}
\end{equation}
Similar to the argument of theorem~\ref{39}, it can be seen that the last line of equation (2.6) goes to 0. Hence
$\|(T_{f}-T_{\delta_{\varepsilon}\ast f_{r_{n}}})g\|_{2,\alpha}\rightarrow 0$ as $r_{n}\rightarrow1$, $\varepsilon\rightarrow0$ and for all $\varepsilon_{1}>0$. This finishes the demonstration that $T_{\delta_{\varepsilon}\ast f_{r_{n}}}\rightarrow T_{f}$ in the strong operator topology.

 The reader may show that $T_{\delta_{\varepsilon}\ast f_{r_{n}}}^{*}\rightarrow T_{f}^{*}$ in the strong operator topology by an argument like the ones previously given.

From Theorem~\ref{40}, there are two sequences $\{a_{n}\}_{n\geq1}$ and $\{b_{n}\}_{n\geq1}$ of non-negative real numbers such that $\sum_{n\geq1}a_{n}=\sum_{n\geq1}b_{n}=1$. Set $\varphi_{1}=\sum^{\infty}_{n=1}a_{n}\psi_{n}$ and  $\varphi_{2}=\sum^{\infty}_{n=1}b_{n}\psi_{n}$. Since every $\psi_{n}\in C_{\partial}(\overline{D})$, we have $\varphi_{1},\varphi_{2}\in C_{\partial}(\overline{D})$.

 Hence $T_{\varphi_{1}}$ and $T_{\varphi_{2}}$ are compact and they satisfy $\|T_{f}-T_{\varphi_{1}}\|=\|T_{f}-T_{\varphi_{2}}\|=\|T_{g}\|_{e}$. Theorem~\ref{40} guarantees $T_{\varphi_{1}}\neq T_{\varphi_{2}}$. By setting $\varphi=s\varphi_{1}+(1-s)\varphi_{2}$, for $s\in (0,1)$, we produce
 infinitely many compact Toeplitz operators $T_{\varphi}$ such that $\|T_{f}-T_{\varphi}\|=\|T_{f}\|_{e}$. \qed

 \textbf{ Proof of Theorem~\ref{18}}. Let $\hbar(D)$ denote the collection of all bounded harmonic functions on the unit disk $D$. For $f\in \hbar(D)$, we'll show that there exists a sequence $\{f_{n}\} $ of functions harmonic on $D$ and continuous on $\overline{D}$ such that the compact Hankel operators $H_{f_{n}}$($H_{f_{n}}^{*}$) converge to $H_{f}$($H_{f}^{*}$) in the strong operator topology.

For $f\in \hbar(D)$ and $z\in D$, let
   $$f_{n}(z)=f(r_{n}z)$$
   where $0<r_{n}<1$ and $r_{n}\rightarrow 1$ as $n\rightarrow \infty$. We have $\triangle f_{n}(z)=r_{n}^{2}\triangle f(r_{n}z)=0$ where $\triangle$ is the Laplace operator; thus $f_{n}$ is harmonic on $D$ and continuous on $\overline{D}$, and by \cite[Sec 8.4]{Zhu}, $H_{f_{n}}$ is compact on the weighted Bergman space. It's not hard to see that $f_{n}(z)$ converges to $f(z)$ pointwise as $n\rightarrow \infty$. By the maximum modulus principle, we have $\|f_{n}\|_{\infty}\leq\|f\|_{\infty}$. For $g\in A^{2}_{\alpha}(D)$, $|f_{n}(z)g(z)|\leq\|f\|_{\infty}|g(z)|$ and $f_{n}(z)g(z)$ pointwise converges to $f(z)g(z)$. We can now apply the dominated convergence theorem to see that
   $$\lim_{n\rightarrow \infty}\int_{D}|(f_{n}(z)-f(z))g(z)|^{2}dA_{\alpha}(z)=0.$$
This yields
\begin{equation}
  \begin{aligned}
&\|(H_{f_{n}}-H_{f})g\|_{2,\alpha}^{2}=\|(I-P)(f_{n}-f)g\|_{2,\alpha}^{2}\\
&\leq\|(f_{n}-f)g\|^{2}_{2,\alpha}=\int_{D}|(f_{n}(w)-f(w))g(w)|^{2}dA_{\alpha}(w)\rightarrow 0
 \end{aligned}
\end{equation}
which implies  $H_{f_{n}}\rightarrow H_{f}$ in the strong operator topology.

Similarly, $H_{f_{n}}^{*}\rightarrow H_{f}^{*}$ in the strong operator topology.

Let  $\{a_{n}\}_{n\geq1}$ and $\{b_{n}\}_{n\geq1}$ be two non-negative real valued sequences  such that $\sum_{n\geq1}a_{n}=\sum_{n\geq1}b_{n}=1$, and set $\varphi_{1}=\sum^{\infty}_{n=1}a_{n}f_{n}$ and  $\varphi_{2}=\sum^{\infty}_{n=1}b_{n}f_{n}$. Since each $f_{n}$ is harmonic on $D$ and continuous on $\overline{D}$, we have $\varphi_{1}$ and $\varphi_{2}$ are  harmonic on $D$ and continuous on $\overline{D}$;
 thus, $H_{\varphi_{1}}$ and $H_{\varphi_{2}}$ are compact operators which satisfy $\|H_{f}-H_{\varphi_{1}}\|=\|H_{f}-H_{\varphi_{2}}\|=\|H_{f}\|_{e}$, and are distinct by Theorem~\ref{40}. The proof can be finished in the now standard fashion.\qed

 $\mathbf{Acknowledgement}$: I would like to thank Prof. Dechao Zheng for introducing the subject of Hankel and Toplitz operators to me, and for his guidance in all aspects of this project. This research was conducted while I was a visiting student at Vanderbilt University. I would like to thank Vanderbilt University for their hospitality.\\

\bigskip

\noindent
{Fengying Li} \hfill \textsl{e-mail}: \texttt{lify0308@163.com}\\
{\scshape
Department of Mathematics \\
Sichuan University \\
Chengdu, Sichuan 610064, China}

\end{document}